\documentclass[11pt]{article}

\usepackage{latexsym}
\usepackage{amssymb}

\newtheorem{thm}{Theorem}
\newtheorem{prop}{Proposition}

\newtheorem{lem}{Lemma}
\newtheorem{dfn}{Definition}

\begin{document}
{
\begin{center}
{\Large\bf
Characteristic properties for a generalized resolvent of a pair of commuting isometric operators.}
\end{center}
\begin{center}
{\bf S.M. Zagorodnyuk}
\end{center}

\section{Introduction.}

Let $V_1, V_2$ be closed isometric operators in a Hilbert space $H$. Suppose that
\begin{equation}
\label{f1_1}
V_1 V_2 h = V_2 V_1 h,\quad h\in D(V_1 V_2)\cap D(V_2 V_1).
\end{equation}
In general, it is not an easy question whether there exist a Hilbert space $\widetilde H\supseteq H$ and
commuting unitary operators $U_1,U_2$ in $\widetilde H$, such that $U_1\supseteq V_1$, $U_2\supseteq V_2$.
This problem was studied in a series of papers \cite{cit_4000_Koranyi}, \cite{cit_6000_Markelov}, \cite{cit_700_Arocena},
\cite{cit_6500_Moran}, \cite{cit_7000_Moran}, see also references therein.
If the answer on the above question is affirmative, then we may define the following operator-valued function of
two complex variables:
$$ \mathbf{R}_{z_1,z_2} = \mathbf{R}_{z_1,z_2} (V_1,V_2) = $$
$$ =
\left. 
P^{\widetilde H}_H (E_{\widetilde H} + z_1 U_1) (E_{\widetilde H} - z_1 U_1)^{-1}
(E_{\widetilde H} + z_2 U_2) (E_{\widetilde H} - z_2 U_2)^{-1} \right|_H, $$
\begin{equation}
\label{f1_5}
z_1, z_2\in \mathbb{T}_e.
\end{equation}
The function $\mathbf{R}_{z_1,z_2} (V_1,V_2)$ is called \textit{a generalized resolvent of a pair of isometric operators
$V_1, V_2$ (corresponding to extensions $U_1,U_2$)}.
Let $\widetilde E_{k,t}$, $t\in[0,2\pi]$, be the (right-continuous) spectral family~\footnote{We shall use the terminology
from~\cite{cit_7500_R_S_N}.} of $U_k$, 
$\widetilde E_{k,0} =0$,
$k=1,2$. 
The following
operator-valued function of two real variables:
\begin{equation}
\label{f1_7}
\mathbf{E}_{t_1,t_2} =
\left. 
P^{\widetilde H}_H \widetilde E_{1,t_1} \widetilde E_{2,t_2} \right|_H,\qquad t_1,t_2\in [0,2\pi],
\end{equation}
is said to be \textit{a (strongly right-continuous) spectral function of a pair of isometric operators
$V_1, V_2$ (corresponding to extensions $U_1,U_2$)}.
As it follows from their definitions, a generalized resolvent and a spectral function, which correspond to the same
extensions $U_1,U_2$, are related by the following equality:
$$ \left(
\mathbf{R}_{z_1,z_2} h,h
\right)_H
=
\int_{\mathbb{R}^2}
\left(
\frac{1+z_1 e^{it_1}}{1-z_1 e^{it_1}}
\right)
\left(
\frac{1+z_2 e^{it_2}}{1-z_2 e^{it_2}}
\right)
d(
\mathbf{E}_{t_1,t_2} h,h )_H, $$
\begin{equation}
\label{f1_9}
h\in H,\ z_1,z_2\in \mathbb{T}_e.
\end{equation}
Here the "distribution" function $(
\mathbf{E}_{t_1,t_2} h,h )_H$ 
defines a (non-negative) finite measure $\sigma$ on $\mathfrak{B}(\mathbb{R}^2)$.
Moreover, we have $\sigma((0,2\pi]\times(0,2\pi]) = \sigma(\mathbb{R}^2) = \| h \|^2_H$.
(One may define $\sigma$ on a semi-ring of rectangles of the form
$\delta = \{ a < t_1 \leq b,\ c < t_2 \leq d \}$
and then extend by the standard procedure).


Let $V$ be a closed isometric operator in a Hilbert space $H$. Then there always exists a unitary extension $U\supseteq V$
in a Hilbert space $\widetilde H\supseteq H$.
Recall that the following operator-valued function:
\begin{equation}
\label{f1_12}
\mathbf{R}_\zeta = \mathbf{R}_\zeta(V) = 
\left. P^{\widetilde{H}}_H \left( E_{\widetilde{H}} - \zeta U \right)^{-1} \right|_H,\qquad
\zeta\in \mathbb{T}_e,
\end{equation}
is said to be  \textit{a generalized resolvent
of an isometric operator $V$ (corresponding to the extension $U$)}.
An arbitrary generalized resolvent $\mathbf{R}_\zeta$ has the following form (\cite{cit_750_Ch}):
\begin{equation}
\label{f1_15}
\mathbf{R}_\zeta = \left[ E_H - \zeta (V\oplus F_\zeta) \right]^{-1},\qquad \zeta\in \mathbb{D},
\end{equation}
where $F_\zeta$ is a function from $\mathcal{S}(\mathbb{D};N_0(V),N_\infty(V))$.
Conversely, an arbitrary function $F_\zeta\in\mathcal{S}(\mathbb{D};N_0(V),N_\infty(V))$
defines by relation~(\ref{f1_15}) a generalized resolvent $\mathbf{R}_\zeta$ of the operator $V$.
Moreover, to different functions from $\mathcal{S}(\mathbb{D};N_0(V),N_\infty(V))$ there correspond
different generalized resolvents of the operator $V$.
Formula~(\ref{f1_15}) is known as Chumakin's formula for the generalized resolvents of an isometric operator.
Moreover, Chumakin established the following characteristic properties of a generalized resolvent
of a closed isometric operator (\cite{cit_750_Ch}): 

\begin{thm}
\label{t1_2_p1_1}

In order that a family of linear operators $R_\zeta$, acting in a Hilbert space $H$ ($D_{R_\zeta} = H$)
and depending on complex parameter $\zeta$ ($|\zeta|\not= 1$), be a generalized resolvent of a closed isometric operator,
it is necessary and sufficient that the following conditions hold:

\begin{itemize}

\item[1)] There exists a number $\zeta_0\in \mathbb{D}\backslash\{ 0 \}$ and a subspace
$L\subseteq H$ such that
$$ ( \zeta R_\zeta - \zeta_0 R_{\zeta_0} ) f
= (\zeta-\zeta_0) R_\zeta R_{\zeta_0} f, $$
for arbitrary $\zeta\in \mathbb{T}_e$ and $f\in L$;

\item[2)] The operator $R_0$ is bounded and $R_0 h = h$, for all $h\in H\ominus \overline{ R_{\zeta_0}L }$;

\item[3)] For an arbitrary $h\in H$ the following inequality holds:
$$ \mathop{\rm Re}\nolimits (R_\zeta h,h)_H \geq \frac{1}{2} \| h \|_H^2,\qquad \zeta\in \mathbb{D}; $$

\item[4)] For an arbitrary $h\in H$
$R_\zeta h$ is an analytic vector-valued function of a parameter $\zeta$
in $\mathbb{D}$;

\item[5)] For an arbitrary $\zeta\in \mathbb{D}\backslash\{ 0 \}$ holds:
$$ R_\zeta^* = E_H - R_{\frac{1}{ \overline{\zeta} }}. $$
\end{itemize}
\end{thm}

\begin{thm}
\label{t1_3_p1_1}

In order that a family of linear operators $R_\zeta$ ($D_{R_\zeta} = H$, $|\zeta|\not= 1$) in a Hilbert space $H$ be a generalized resolvent 
of a given closed isometric operator $V$ in $H$,
it is necessary and sufficient that the following conditions hold:


\begin{itemize}

\item[1)] For all $\zeta\in \mathbb{T}_e$ and for all $g\in D(V)$ the following equality holds:
$$ R_\zeta (E_H - \zeta V) g = g; $$

\item[2)] The operator $R_0$ is bounded and $R_0 h = h$, for all $h\in H\ominus D(V)$;

\item[3)] For an arbitrary $h\in H$ the following inequality holds:
$$ \mathop{\rm Re}\nolimits (R_\zeta h,h)_H \geq \frac{1}{2} \| h \|_H^2,\qquad \zeta\in \mathbb{D}; $$

\item[4)] For an arbitrary $h\in H$
$R_\zeta h$ is an analytic vector-valued function of a parameter $\zeta$
in $\mathbb{D}$;

\item[5)] For an arbitrary  $\zeta\in \mathbb{D}\backslash\{ 0 \}$ the following equality is true:
$$ R_\zeta^* = E_H - R_{\frac{1}{ \overline{\zeta} }}. $$
\end{itemize}
\end{thm}

Our purpose is to obtain an analog of Theorem~\ref{t1_2_p1_1} for a generalized resolvent of
a pair of commuting isometric operators. An important role will be played by the following class $H_2$ of analytic functions of
two complex variables, which was introduced by Kor\'anyi in~\cite{cit_4000_Koranyi} (We use the original notation of
Kor\'anyi for this class. Since the Hardy space will not appear in this paper, it will cause no confusion).

\begin{dfn}
\label{d1_1}
The class $H_2$ is the class of functions $f$ of two complex variables $z_1,z_2$ defined and holomorphic for all
$|z_1|,|z_2|\not= 1$ (including $\infty$) and satisfying the conditions
\begin{itemize}

\item[(a)] $f(\overline{z_1}^{-1}, \overline{z_2}^{-1}) = \overline{f}(z_1,z_2)$ for all  $|z_1|,|z_2|\not= 1$,

\item[(b)] $f(z_1,z_2) - f(\overline{z_1}^{-1}, z_2) - f(z_1, \overline{z_2}^{-1}) + f(\overline{z_1}^{-1},\overline{z_2}^{-1}) \geq 0$, 
for $|z_1|,|z_2| < 1$,

\item[(c)] $f(z_1,0) + f(z_1,\infty) = 0$, $f(0,z_2) + f(\infty,z_2) = 0$ for all $|z_1|\not= 1$ and $|z_2|\not= 1$. 

\end{itemize}
\end{dfn}

Every function $g\in H_2$ admits the following representation (see~\cite[formula (26)]{cit_4000_Koranyi} and considerations on page~532 
in~\cite{cit_4000_Koranyi}):
$$ g(z_1,z_2) = \frac{1}{4} \left(
(E+z_1 \widehat U)(E-z_1 \widehat U)^{-1} (E+z_2 \widehat V)(E-z_2 \widehat V)^{-1} \varepsilon_{0,0}, \varepsilon_{0,0}
\right)_{ \widehat{\mathfrak{B}} }, $$
\begin{equation}
\label{f1_18}
z_1,z_2\in \mathbb{T}_e,
\end{equation}
where $\widehat U, \widehat V$ are some commutative unitary operators in a Hilbert space $\widehat{\mathfrak{B}}$; 
$\varepsilon_{0,0}\in \widehat{\mathfrak{B}}$.
Let $\widehat E_{1,t}$, $t\in[0,2\pi]$, be the (right-continuous) spectral family of $\widehat U$, 
$\widehat E_{1,0} =0$. 
Let $\widehat E_{2,t}$, $t\in[0,2\pi]$, be the (right-continuous) spectral family of $\widehat V$, 
$\widehat E_{2,0} =0$.
As in relation~(\ref{f1_9}) we may write:
$$ g(z_1,z_2) = 
\int_{\mathbb{R}^2}
\left(
\frac{1+z_1 e^{it_1}}{1-z_1 e^{it_1}}
\right)
\left(
\frac{1+z_2 e^{it_2}}{1-z_2 e^{it_2}}
\right)
d\left(
\widehat E_{1,t_1} \widehat E_{2,t_2} \frac{1}{2}\varepsilon_{0,0}, \frac{1}{2} \varepsilon_{0,0} \right)_{ \widehat{\mathfrak{B}} };  $$
\begin{equation}
\label{f1_20}
g(z_1,z_2) = 
\int_{\mathbb{R}^2}
\left(
\frac{1+z_1 e^{it_1}}{1-z_1 e^{it_1}}
\right)
\left(
\frac{1+z_2 e^{it_2}}{1-z_2 e^{it_2}}
\right)
d\mu,\quad z_1,z_2\in \mathbb{T}_e,
\end{equation}
where $\mu$ is a (non-negative) finite measure on $\mathfrak{B}(\mathbb{R}^2)$ generated by the distribution function
$\left( \widehat E_{1,t_1} \widehat E_{2,t_2} \frac{1}{2} \varepsilon_{0,0}, \frac{1}{2} \varepsilon_{0,0} \right)_{ \widehat{\mathfrak{B}} }$. 
Moreover, we have $\mu((0,2\pi]\times(0,2\pi]) = \mu(\mathbb{R}^2)$.


Another important ingredient of our proof is generalized Neumark's dilation theorem~\cite[p. 499]{cit_7500_R_S_N} (while in the proof of
Chumakin's result the usual Neumark's dilation theorem is used).

\noindent
\textbf{Notations.}
As usual, we denote by $\mathbb{R}, \mathbb{C}, \mathbb{N}, \mathbb{Z}, \mathbb{Z}_+$,
the sets of real numbers, complex numbers, positive integers, integers and non-negative integers,
respectively; $\mathbb{D} = \{ z\in \mathbb{C}:\ |z|<1 \}$; $\mathbb{D}_e = \{ z\in \mathbb{C}:\ |z|>1 \}$;
$\mathbb{T} = \{ z\in \mathbb{C}:\ |z|=1 \}$; $\mathbb{T}_e = \{ z\in \mathbb{C}:\ |z|\not= 1 \}$.
By $k\in\overline{m,n}$ (or $k = \overline{m,n}$) we mean that $k\in \mathbb{Z}_+:\ m\leq k\leq n$; for $m,n\in\mathbb{Z}_+$.
By $\mathbb{R}^2$ we denote the two-dimensional real Eucledian space.
By $\mathfrak{B}(\mathbb{R}^2)$ we mean the set of all Borel subsets of $\mathbb{R}^2$.

\noindent
In this paper Hilbert spaces are not necessarily separable, operators in them are
supposed to be linear.

\noindent
If H is a Hilbert space then $(\cdot,\cdot)_H$ and $\| \cdot \|_H$ mean
the scalar product and the norm in $H$, respectively.
Indices may be omitted in obvious cases.
For a linear operator $A$ in $H$, we denote by $D(A)$
its  domain, by $R(A)$ its range, and $A^*$ means the adjoint operator
if it exists. If $A$ is invertible then $A^{-1}$ means its
inverse. $\overline{A}$ means the closure of the operator, if the
operator is closable. If $A$ is bounded then $\| A \|$ denotes its
norm.
For a set $M\subseteq H$
we denote by $\overline{M}$ the closure of $M$ in the norm of $H$.
By $\mathop{\rm Lin}\nolimits M$ we mean
the set of all linear combinations of elements from $M$,
and $\mathop{\rm span}\nolimits M:= \overline{ \mathop{\rm Lin}\nolimits M }$.
By $E_H$ we denote the identity operator in $H$, i.e. $E_H x = x$,
$x\in H$. In obvious cases we may omit the index $H$. If $H_1$ is a subspace of $H$, then $P_{H_1} =
P_{H_1}^{H}$ is an operator of the orthogonal projection on $H_1$
in $H$.
By $[H]$ we denote a set of all bounded operators on $H$.
For a closed isometric operator $V$ in $H$ we denote:
$M_\zeta(V) = (E_H - \zeta V) D(V)$, $N_\zeta(V) = H\ominus M_\zeta(V)$, $\zeta\in \mathbb{C}$; $M_\infty(V)=R(V)$, $N_\infty(V)= H\ominus R(V)$.
For a unitary operator $U$ in $H$ we denote:
$\mathcal{R}_{z}(U) := ( E_{H} - z U )^{-1}$, $z\in\mathbb{T}_e$.

By $\mathcal{S}(D;N,N')$ we denote a class of all analytic in a domain $D\subseteq \mathbb{C}$
operator-valued functions $F(z)$, which values are linear non-expanding operators mapping the whole
$N$ into $N'$, where $N$ and $N'$ are some Hilbert spaces.

For a unitary operator $U$ in a Hilbert space $H$ we shall use the following notation:
$$ U(z) := (E_{H} + z U) (E_{H} - z U)^{-1} = -E_{H} + 2 \mathcal{R}_{z}(U),\qquad z\in\mathbb{T}_e. $$
It is straightforward to check that (\cite[p. 531]{cit_4000_Koranyi})
\begin{equation}
\label{f4_5}
U^*(z) = -U\left(
\frac{1}{ \overline{z} }
\right),\qquad  z\in\mathbb{T}_e\backslash\{ 0 \};
\end{equation}
\begin{equation}
\label{f4_7}
U(z) - U\left(
\frac{1}{ \overline{z} }
\right) =
2 (1-|z|^2) \mathcal{R}_{z}^*(U) \mathcal{R}_{z}(U) \geq 0,\qquad z\in\mathbb{D}\backslash\{ 0 \}.
\end{equation}
If we set $U(\infty) := - E_H$, then relation~(\ref{f4_5}) will be valid for all $z\in \mathbb{T}_e\cup\{ \infty \}$.

\section{Preliminary results.}

We shall need the following elementary lemma.

\begin{lem}
\label{l2_1}
Let $\mu$ be a (non-negative) finite measure on $\mathfrak{B}(\mathbb{R}^2)$.
Let $\varphi_j(z;t)$ be an analytic of $z$ in a domain $D\subseteq\mathbb{C}$ complex-valued function
depending on a parameter $t\in\mathbb{R}$
with all derivatives $(\varphi_j(z;t))^{(k)}_{z}$, $k\in\mathbb{Z}_+$ being continuous and bounded as a function of $t$
(with an arbitrary fixed $z\in D$); $j=1,2$.
Suppose that for each $z_0\in D$ there exists a closed ball $U(z_0) = \{ z\in\mathbb{C}:\ |z-z_0|\leq R_{z_0} \} \subseteq D$
($R_{z_0}> 0$), such that  
\begin{equation}
\label{f2_1_0}
\left|
(\varphi_j(z;t))^{(k)}_{z}
\right| \leq M_{k,j}(z_0),\qquad z\in U(z_0),\ t\in\mathbb{R},\ k\in\mathbb{Z}_+,
\end{equation}
where $M_{k,j}(z_0)$ does not depend on $t$. Here $j=1,2$ is a fixed number.
Then
\begin{equation}
\label{f2_1}
\left( \left( g(z_1,z_2) \right)^{(k)}_{z_1} \right)^{(l)}_{z_2} =
\int_{\mathbb{R}^2} \left( \varphi_1(z_1;t_1) \right)^{(k)}_{z_1} \left( \varphi_2(z_2;t_2) \right)^{(l)}_{z_2} d\mu(t_1,t_2),\quad
k,l\in\mathbb{Z}_+,
\end{equation}
where
\begin{equation}
\label{f2_2}
g(z_1,z_2) =
\int_{\mathbb{R}^2} \varphi_1(z_1;t_1) \varphi_2(z_2;t_2) d\mu(t_1,t_2),\quad
z_1,z_2\in D,
\end{equation}
and all derivatives in~(\ref{f2_1}) exist.
\end{lem}   
\textbf{Proof.}
Firstly, we shall check relation~(\ref{f2_1}) with $l=0$ by the induction (for $k\in\mathbb{Z}_+$).
We may use the definition of the derivative, 
Lagrange's theorem on a finite increament of a function (the mean value theorem), inequality~(\ref{f2_1_0}) 
and the Lebesgue dominated convergence theorem
to verify the induction step.
Secondly, fix an arbitrary $k\in\mathbb{Z}_+$ and check relation~(\ref{f2_1}) by the induction (for $l\in\mathbb{Z}_+$)
in a similar manner.
$\Box$

By the induction argument we may write:
\begin{equation}
\label{f2_3}
\left(
\frac{1+z e^{it}}{1-z e^{it}}
\right)^{(k)}_z = 2 k!\frac{ e^{ikt} }{ (1-z e^{it})^{k+1} } - \delta_{k,0},\quad
z\in \mathbb{T}_e,\ t\in \mathbb{R},\ k\in \mathbb{Z}_+;
\end{equation}
\begin{equation}
\label{f2_4}
\left(
\frac{u + e^{it}}{u - e^{it}}
\right)^{(l)}_u = (-1)^l 2 l!\frac{ e^{it} }{ (u-e^{it})^{l+1} } + \delta_{l,0},\quad
u\in \mathbb{T}_e,\  t\in \mathbb{R},\ l\in \mathbb{Z}_+.
\end{equation}

Let $g(z_1,z_2)$ be an arbitrary function which admits representation~(\ref{f1_20})
where $\mu$ is a (non-negative) finite measure on $\mathfrak{B}(\mathbb{R}^2)$ with $\mu((0,2\pi]\times(0,2\pi]) = \mu(\mathbb{R}^2)$.
By Lemma~\ref{l2_1} and relations~(\ref{f2_3}),(\ref{f2_4}) we obtain that
\begin{equation}
\label{f2_7}
\left. \left( \left( g(z_1,z_2) \right)^{(k)}_{z_1} \right)^{(l)}_{z_2} \right|_{(z_1,z_2) = (0,0)} =
\left\{
\begin{array}{cccc} s_{0,0}, & \mbox{if } k=l=0\\
2 l! s_{0,l}, & \mbox{if } k=0,\ l\in\mathbb{N}\\
2 k! s_{k,0}, & \mbox{if } k\in\mathbb{N},\ l=0\\
4 k! l! s_{k,l}, & \mbox{if } k,l\in\mathbb{N}\end{array}
\right.;
\end{equation}
\begin{equation}
\label{f2_9}
\left. \left( \left( g(u_1^{-1},z_2) \right)^{(k)}_{u_1} \right)^{(l)}_{z_2} \right|_{(u_1,z_2) = (0,0)} =
\left\{
\begin{array}{cc} 
- 2 k! s_{-k,0}, & \mbox{if } k\in\mathbb{N},\ l=0\\
-4 k! l! s_{-k,l}, & \mbox{if } k,l\in\mathbb{N} \end{array}
\right.,
\end{equation}
where $g(u_1^{-1},z_2)|_{u_1=0} := \lim_{u_1\to 0} g(u_1^{-1},z_2)$, $z_2\in\mathbb{D}$; and therefore $g(u_1^{-1},z_2)$ is defined on
$\mathbb{D}\times\mathbb{D}$;
\begin{equation}
\label{f2_12}
\left. \left( \left( g(z_1,u_2^{-1}) \right)^{(k)}_{z_1} \right)^{(l)}_{u_2} \right|_{(z_1,u_2) = (0,0)} =
\left\{
\begin{array}{cccc} 
- 2 l! s_{0,-l}, & \mbox{if } k=0,\ l\in\mathbb{N}\\
- 4 k! l! s_{k,-l}, & \mbox{if } k,l\in\mathbb{N}\end{array}
\right.,
\end{equation}
where $g(z_1,u_2^{-1})|_{u_2=0} := \lim_{u_2\to 0} g(z_1,u_2^{-1})$, $z_1\in\mathbb{D}$; and therefore $g(z_1,u_2^{-1})$ is defined on
$\mathbb{D}\times\mathbb{D}$;
\begin{equation}
\label{f2_17}
\left. \left( \left( g(u_1^{-1},u_2^{-1}) \right)^{(k)}_{u_1} \right)^{(l)}_{u_2} \right|_{(u_1,u_2) = (0,0)} =
4 k! l! s_{-k,-l},\quad k,l\in\mathbb{N},
\end{equation}
where $g(u_1^{-1},u_2^{-1})|_{u_1=0} = \lim_{u_1\to 0} g(u_1^{-1},u_2^{-1})$, $u_2\in\mathbb{D}\backslash\{ 0 \}$;
$g(u_1^{-1},u_2^{-1})|_{u_2=0} = \lim_{u_2\to 0} g(u_1^{-1},u_2^{-1})$, $u_1\in\mathbb{D}\backslash\{ 0 \}$;
$g(u_1^{-1},u_2^{-1})|_{u_1=u_2=0} = \lim_{u_2\to 0} g(u_1^{-1},u_2^{-1})|_{u_1=0}$; and therefore $g(u_1^{-1},u_2^{-1})$ is defined on
$\mathbb{D}\times\mathbb{D}$.
Here
\begin{equation}
\label{f2_22}
s_{k,l} := 
\int_{\mathbb{R}^2}
e^{ikt_1}
e^{ilt_2}
d\mu,\quad k,l\in \mathbb{Z},
\end{equation}
are the trigonometric moments of $\mu$.
Thus, \textit{all trigonometric moments of $\mu$ are uniquely determined by the function $g(z_1,z_2)$}.

Consider the following function:
\begin{equation}
\label{f2_23}
f_{m,k}(t) =
\left\{
\begin{array}{cc}
\left(
\left( \frac{1}{k} \right)^m - (2\pi)^m
\right)
kt + (2\pi)^m, & 0\leq t\leq \frac{1}{k} \\
t^m, & \frac{1}{k} < t \leq 2\pi
\end{array}
\right.,
\end{equation}
where $m\in\mathbb{Z}_+$, $k\in\mathbb{N}$.
Extend $f_{m,k}(t)$ to a continuous function on the real line with the period $2\pi$.
By Weierstrass's approximation theorem there exists a trigonometric polynomial $T_{m,k}(t)$ such that
\begin{equation}
\label{f2_23_1}
\left| f_{m,k}(t) - T_{m,k}(t) \right| < \frac{1}{k},\qquad t\in\mathbb{R}.
\end{equation}
Observe that 
\begin{equation}
\label{f2_23_2}
\left| f_{m,k}(t) \right| \leq (2\pi)^m,\qquad t\in\mathbb{R}.
\end{equation}
By~(\ref{f2_23_1}) it follows that
\begin{equation}
\label{f2_23_3}
\left| T_{m,k}(t) \right| \leq (2\pi)^m +1,\qquad t\in\mathbb{R}.
\end{equation}
For arbitrary $m,n\in\mathbb{Z}_+$ we may write
$$ \left|
\int_{\mathbb{R}^2} t_1^m t_2^n d\mu - \int_{\mathbb{R}^2} T_{m,k}(t_1) T_{n,k}(t_2) d\mu
\right| \leq $$
$$ \leq 
\left|
\int_{\mathbb{R}^2} \left( t_1^m - T_{m,k}(t_1) \right) t_2^n d\mu \right| + 
\left|
\int_{\mathbb{R}^2} T_{m,k}(t_1) \left( t_2^n - T_{n,k}(t_2) \right) d\mu
\right| \leq $$
$$ \leq 
\left|
\int_{\mathbb{R}^2} \left( t_1^m - f_{m,k}(t_1) \right) t_2^n d\mu \right| + 
\left|
\int_{\mathbb{R}^2} \left( f_{m,k}(t_1) - T_{m,k}(t_1) \right) t_2^n d\mu
\right| + 
$$
$$ + \left|
\int_{\mathbb{R}^2} T_{m,k}(t_1) \left( t_2^n - f_{n,k}(t_2) \right) d\mu \right| 
+ $$
\begin{equation}
\label{f2_24}
+ \left|
\int_{\mathbb{R}^2} T_{m,k}(t_1) \left( f_{n,k}(t_2) - T_{n,k}(t_2) \right) d\mu \right| 
\rightarrow 0, 
\end{equation}
as $k\rightarrow\infty$.
Therefore all power moments:
\begin{equation}
\label{f2_25}
r_{m,n} := 
\int_{\mathbb{R}^2}
t_1^m t_2^n
d\mu,\quad m,n\in \mathbb{Z_+},
\end{equation}
are uniquely determined by the function $g(z_1,z_2)$.
Since the two-dimensional power moment problem which has a solution with a compact support is determinate
(e.g.~\cite[Theorem B, p. 323]{cit_6020_McGregor}), then
we conclude that \textit{the measure $\mu$ in representation~(\ref{f1_20}) is uniquely determined by
the function $g$}.

\begin{prop}
\label{p2_1}
Let $\sigma_j$ ($j=\overline{1,4}$) be (non-negative) finite measures on $\mathfrak{B}(\mathbb{R}^2)$ with
$\sigma_j((0,2\pi]^2) = \sigma_j(\mathbb{R}^2)$.
If
\begin{equation}
\label{f2_30}
s_{k,l}(\sigma_1) - s_{k,l}(\sigma_2) + i s_{k,l}(\sigma_3) - i s_{k,l}(\sigma_4) = 0,\qquad
k,l\in\mathbb{Z},
\end{equation}
then
\begin{equation}
\label{f2_32}
\sigma_1 - \sigma_2 + i \sigma_3 - i \sigma_4 = 0.
\end{equation}
\end{prop}
\textbf{Proof.}
Observe that the measures $\sigma_j$ ($j=\overline{1,4}$) satisfy the assumptions on the measure $\mu$ introduced after~(\ref{f2_4}).
Therefore we may apply the above constructions to these measures. Notice that the function $f_{m,k}(t)$ in~(\ref{f2_23})
depends on $m,k,t$ but do not depend on the measure $\mu$.
By~(\ref{f2_24}) for arbitrary $m,n\in\mathbb{Z}_+$ we may write
$$ \left|
r_{m,n}(\sigma_1) - r_{m,n}(\sigma_2) + i r_{m,n}(\sigma_3) - i r_{m,n}(\sigma_4) -
\left(
\int_{\mathbb{R}^2} T_{m,k}(t_1) T_{n,k}(t_2) d\sigma_1 -
\right.
\right.
$$
$$  
- \int_{\mathbb{R}^2} T_{m,k}(t_1) T_{n,k}(t_2) d\sigma_2 +
i \int_{\mathbb{R}^2} T_{m,k}(t_1) T_{n,k}(t_2) d\sigma_3 -
$$
\begin{equation}
\label{f2_35}
- \left. \left.
i \int_{\mathbb{R}^2} T_{m,k}(t_1) T_{n,k}(t_2) d\sigma_4
\right)
\right| \rightarrow 0, 
\end{equation}
as $k\rightarrow\infty$.
By~(\ref{f2_30}) we conclude that the expression in the round brackets in~(\ref{f2_35}) is equal to zero.
Therefore
\begin{equation}
\label{f2_37}
r_{m,n}(\sigma_1) - r_{m,n}(\sigma_2) + i r_{m,n}(\sigma_3) - i r_{m,n}(\sigma_4) = 0,\qquad
m,n\in\mathbb{Z}_+.
\end{equation}
Extracting the real and the imaginary parts we get
\begin{equation}
\label{f2_39}
r_{m,n}(\sigma_1) = r_{m,n}(\sigma_2),\qquad m,n\in\mathbb{Z}_+;
\end{equation}
\begin{equation}
\label{f2_41}
r_{m,n}(\sigma_3) = r_{m,n}(\sigma_4),\qquad m,n\in\mathbb{Z}_+.
\end{equation}
Since the corresponding two-dimensional power moment problem is determinate, we conclude that $\sigma_1 = \sigma_2$
and $\sigma_3 = \sigma_4$.
$\Box$

\begin{prop}
\label{p2_2}
Let $\sigma_j$ ($j=\overline{1,4}$) be (non-negative) finite measures on $\mathfrak{B}(\mathbb{R}^2)$ with
$\sigma_j((0,2\pi]^2) = \sigma_j(\mathbb{R}^2)$.
Let $g_j(z_1,z_2)$ be a function which admits representation~(\ref{f1_20}) with $\sigma_j$ instead of $\mu$;
$j=\overline{1,4}$.
If
\begin{equation}
\label{f2_50}
g_1(z_1,z_2) - g_2(z_1,z_2) + i g_3(z_1,z_2) - i g_4(z_1,z_2) = 0,\qquad z_1,z_2\in\mathbb{T}_e,
\end{equation}
then
\begin{equation}
\label{f2_52}
\sigma_1 - \sigma_2 + i \sigma_3 - i \sigma_4 = 0.
\end{equation}
\end{prop}
\textbf{Proof.}
The measures $\sigma_j$ ($j=\overline{1,4}$) satisfy the assumptions on the measure $\mu$ introduced after~(\ref{f2_4}).
Moreover, the functions $g_j(z_1,z_2)$ for $\sigma_j$ are introduced in the same way as $g(z_1,z_2)$ for $\mu$.
Calculating derivatives of
$g_1(z_1,z_2) - g_2(z_1,z_2) + i g_3(z_1,z_2) - i g_4(z_1,z_2)$ 
at various points and
using relations~(\ref{f2_7})-(\ref{f2_17}) we obtain that
$$ s_{k,l}(\sigma_1) - s_{k,l}(\sigma_2) + i s_{k,l}(\sigma_3) - i s_{k,l}(\sigma_4) = 0,\qquad
k,l\in\mathbb{Z}. $$
By Proposition~\ref{p2_1} we conclude that relation~(\ref{f2_52}) holds.
$\Box$

\section{Properties of generalized resolvents.}

The following theorem is an analog of Theorem~\ref{t1_2_p1_1}.

\begin{thm}
\label{t3_1}
Let an operator-valued function  $R_{z_1,z_2}$ be given, which depends on complex parameters
$z_1,z_2\in \mathbb{T}_e$ and which values are linear bounded operators defined on a (whole) Hilbert space $H$.
This function is a generalized resolvent of a pair of
closed isometric operators in $H$ (satisfying the commutativity relation~(\ref{f1_1})) if an only if the following conditions are satisfied:

\begin{itemize}

\item[1)] $R_{0,0} = E_H$;

\item[2)] $R_{z_1,z_2}^* = R_{ \frac{1}{\overline{z_1}}, \frac{1}{\overline{z_2}} }$, $z_1,z_2\in \mathbb{T}_e\backslash\{ 0 \}$;

\item[3)] For all $h\in H$, for the function $f(z_1,z_2) := (R_{z_1,z_2} h,h)_H$, $z_1,z_2\in\mathbb{T}_e$, there exist limits:
$$ f(\infty, z_2) := \lim_{z_1\to\infty} f(z_1,z_2),\quad f(z_1,\infty) := \lim_{z_2\to\infty} f(z_1,z_2),\quad z_1,z_2\in \mathbb{T}_e; $$
$$ f(\infty,\infty) = \lim_{z_2\to\infty} \lim_{z_1\to\infty} f(z_1,z_2), $$
and the extended by these relations function $f(z_1,z_2)$, $z_1,z_2\in \mathbb{T}_e\cup\{ \infty \}$ belongs to $H_2$.

\end{itemize}
\end{thm}

\textbf{Proof.} \textit{Necessity.}
Let $V_1, V_2$ be closed isometric operators in a Hilbert space $H$ satisfying relation~(\ref{f1_1}). Suppose that
there exist commuting unitary extensions $U_k\supseteq V_k$, $k=1,2$, in a Hilbert space $\widetilde H\supseteq H$, and
$R_{z_1,z_2} = \mathbf{R}_{z_1,z_2}$ be the corresponding generalized resolvent. 
By the definition of the generalized resolvent we see that condition~1) is satisfied.
By~(\ref{f4_5}) for arbitrary $z_1,z_2\in \mathbb{T}_e\backslash\{ 0 \}$ and $h,g\in H$ we may write
$$ (\mathbf{R}_{z_1,z_2} h, g)_H = \left(
P^{\widetilde H}_H U_1(z_1) U_2(z_2)|_H h, g
\right)_H
=
(U_1(z_1) U_2(z_2) h, g)_{\widetilde H} = $$
$$ = \left( h,
U_1(\overline{z_1}^{-1}) U_2(\overline{z_2}^{-1}) g
\right)_{\widetilde H} = 
(h, \mathbf{R}_{ \frac{1}{\overline{z_1}}, \frac{1}{\overline{z_2}}} g)_H. $$
Therefore condition~2) holds.

\noindent
Choose an arbitrary $h\in H$ and set
\begin{equation}
\label{f3_10}
f(z_1,z_2) = (U_1(z_1) U_2(z_2) h,h)_{\widetilde H},\qquad z_1,z_2\in\mathbb{T}_e\cup\{ \infty \}.
\end{equation}
Here $U_1(\infty) = U_2(\infty) := - E_{\widetilde H}$.
It is easy to check that this definition is consistent with the definition of $f(z_1,z_2)$ from the statement of the
theorem.
Observe that the set $\mathbb{T}_e\times\mathbb{T}_e$ is a union of four polycircular domains $\mathbb{D}\times\mathbb{D}$,
$\mathbb{D}\times\mathbb{D}_e$, $\mathbb{D}_e\times\mathbb{D}$ and $\mathbb{D}_e\times\mathbb{D}_e$.
In each of these domains the function $f(z_1,z_2)$ is holomorphic with respect to each variable.
By Hartogs's theorem we conclude that $f(z_1,z_2)$ is holomorphic at each point of $\mathbb{T}_e\times\mathbb{T}_e$.
For the infinite points we may use the change of variable $u=\frac{1}{z}$ and proceed in the same manner.
Conditions~(a)-(c) in the definition of the class $H_2$ can be checked by relations~(\ref{f4_5}),(\ref{f4_7}), as it was
done in~\cite[p. 531]{cit_4000_Koranyi}. Thus, $f(z_1,z_2)\in H_2$ and condition~3) holds.

\noindent
\textit{Sufficiency.}
Suppose that an operator-valued function $R_{z_1,z_2}$ satisfies the assumptions of the theorem and conditions~1),2),3).
By condition~3) and relation~(\ref{f1_20}) we may write:
$$ ( R_{z_1,z_2} h, h )_H = 
\int_{\mathbb{R}^2}
\left(
\frac{1+z_1 e^{it_1}}{1-z_1 e^{it_1}}
\right)
\left(
\frac{1+z_2 e^{it_2}}{1-z_2 e^{it_2}}
\right)
d\mu(\delta; h,h), $$
\begin{equation}
\label{f3_20}
z_1,z_2\in \mathbb{T}_e,\ h\in H,
\end{equation}
where $\mu(\delta; h,h)$ is a (non-negative) finite measure on $\mathfrak{B}(\mathbb{R}^2)$ such that
$\mu((0,2\pi]\times(0,2\pi]) = \mu(\mathbb{R}^2)$.
Set
$$ \mu(\delta; h,g) = \frac{1}{4} (
\mu(\delta; h+g,h+g) - \mu(\delta; h-g,h-g) + i \mu(\delta; h+ig,h+ig) - $$
\begin{equation}
\label{f3_22}
- i \mu(\delta; h-ig,h-ig) ),\qquad \delta\in \mathfrak{B}(\mathbb{R}^2),\ h,g\in H. 
\end{equation}
Then
$$ ( R_{z_1,z_2} h, g )_H = 
\int_{\mathbb{R}^2}
\left(
\frac{1+z_1 e^{it_1}}{1-z_1 e^{it_1}}
\right)
\left(
\frac{1+z_2 e^{it_2}}{1-z_2 e^{it_2}}
\right)
d\mu(\delta; h,g), $$
\begin{equation}
\label{f3_25}
z_1,z_2\in \mathbb{T}_e,\ h,g\in H.
\end{equation}
The integral of the form $\int_{\mathbb{R}^2} u(t_1,t_2) d\mu(\delta)$ (where $u(t_1,t_2)$ is a complex-valued function on $\mathbb{R}^2$ and
$\mu(\delta)$ is a complex-valued function on $\mathfrak{B}(\mathbb{R}^2)$) may be understood as a limit of Riemann-Stieltjes type integral sums,
if it exists.
This means that we consider partitions of $\mathbb{R}^2$ by rectangles of the following form:
$$ \delta_{n,k} := \{ t_{1,n-1} < t_1 \leq t_{1,n},\ t_{2,k-1} < t_2 \leq t_{2,k} \},\qquad n,k\in\mathbb{Z}, $$
and choose arbitrary points $(t_{1;n,k}, t_{2;n,k})\in \delta_{n,k}$.
The integral sum is defined by $\sum_{n,k} u(t_{1;n,k}, t_{2;n,k}) \mu(\delta_{n,k})$.
The integral is a limit of integral sums as partitions become arbitrarily fine (i.e. the diameter of partitions tends to zero),
if the limit exists,
cf.~\cite[p. 307]{cit_7500_R_S_N}.
 
Fix arbitrary $h,g\in H$.
From the definition of $\mu(\delta; h,g)$ it follows that 
$\mu(\delta; g,h) - \overline{\mu(\delta; h,g)} = \sum_{j=1}^8 \alpha_j \mu_j(\delta)$, $\delta\in \mathfrak{B}(\mathbb{R}^2)$, where
$\alpha_j\in\mathbb{C}$ and $\mu_j(\delta)$ are (non-negative) finite measures on $\mathfrak{B}(\mathbb{R}^2)$
such that $\mu_j((0,2\pi]\times(0,2\pi]) = \mu_j(\mathbb{R}^2)$, $j\in\overline{1,8}$.
Namely, $\{ \alpha_j \}_{j=1}^8 = \{ \frac{1}{4}, -\frac{1}{4}, \frac{i}{4}, -\frac{i}{4},
-\frac{1}{4}, \frac{1}{4}, \frac{i}{4}, -\frac{i}{4} \}$,
$\{ \mu_j \}_{j=1}^8 = \{ \mu(\delta;g+h,g+h), \mu(\delta;g-h,g-h), \mu(\delta;g+ih,g+ih), \mu(\delta;g-ih,g-ih),
\mu(\delta;h+g,h+g), \mu(\delta;h-g,h-g), \mu(\delta;h+ig,h+ig), \mu(\delta;h-ig,h-ig) \}$.
Observe that
$$ \mu_1 = \mu_5,\ \alpha_1 = - \alpha_5;\ 
\mu_2 = \mu_6,\ \alpha_2 = - \alpha_6;\  
\mu_3 = \mu_8,\ \alpha_3 = - \alpha_8; $$
$$ \mu_4 = \mu_7,\ \alpha_4 = - \alpha_7. $$
This follows from the representation~(\ref{f3_20}) for each measure and the established in the previous section fact that
the measure is uniquely determined from the representation of type~(\ref{f1_20}). For example,
$$ (R_{z_1,z_2} (g-ih), g-ih) = (R_{z_1,z_2} (h+ig), h+ig),\qquad z_1,z_2\in\mathbb{T}_e, $$
and therefore $\mu_4=\mu_7$.
Consequently, we obtain the following relation:
\begin{equation}
\label{f3_38}
\mu(\delta; g,h) = \overline{\mu(\delta; h,g)},\quad \delta\in \mathfrak{B}(\mathbb{R}^2),\ h,g\in H.
\end{equation}  
Choose arbitrary $\alpha,\beta\in\mathbb{C}$ and $h_1,h_2,g\in H$. By~(\ref{f3_25}) we may write
$$ \int_{\mathbb{R}^2}
\left(
\frac{1+z_1 e^{it_1}}{1-z_1 e^{it_1}}
\right)
\left(
\frac{1+z_2 e^{it_2}}{1-z_2 e^{it_2}}
\right)
d\mu(\delta; \alpha h_1 + \beta h_2,g) = $$
$$ = ( R_{z_1,z_2} (\alpha h_1 + \beta h_2), g )_H = \alpha ( R_{z_1,z_2} h_1, g )_H + \beta ( R_{z_1,z_2} h_2, g )_H = $$
$$ = 
\alpha
\int_{\mathbb{R}^2}
\left(
\frac{1+z_1 e^{it_1}}{1-z_1 e^{it_1}}
\right)
\left(
\frac{1+z_2 e^{it_2}}{1-z_2 e^{it_2}}
\right)
d \mu(\delta; h_1,g) + $$
$$ +
\beta
\int_{\mathbb{R}^2}
\left(
\frac{1+z_1 e^{it_1}}{1-z_1 e^{it_1}}
\right)
\left(
\frac{1+z_2 e^{it_2}}{1-z_2 e^{it_2}}
\right)
d\mu(\delta; h_2,g),\quad 
z_1,z_2\in\mathbb{T}_e. $$
Therefore
$$ \int_{\mathbb{R}^2}
\left(
\frac{1+z_1 e^{it_1}}{1-z_1 e^{it_1}}
\right)
\left(
\frac{1+z_2 e^{it_2}}{1-z_2 e^{it_2}}
\right)
d( \alpha \mu(\delta; h_1,g) + \beta \mu(\delta; h_2,g) - $$
$$ - \mu(\delta; \alpha h_1 + \beta h_2,g)) = 0,\quad z_1,z_2\in\mathbb{T}_e. $$
By~Proposition~\ref{p2_2} we obtain that
$$ \mu(\delta; \alpha h_1 + \beta h_2,g) = \alpha \mu(\delta; h_1,g) + \beta \mu(\delta; h_2,g), $$
\begin{equation}
\label{f3_42}
\quad \delta\in \mathfrak{B}(\mathbb{R}^2),\ \alpha,\beta\in\mathbb{C},\ 
h_1,h_2,g\in H.
\end{equation}
Observe that 
$$ \left|
\mu(\delta; h,h)
\right| \leq
\mu(\mathbb{R}^2; h,h)
= \int_{ \mathbb{R}^2 }
d\mu(\delta; h,h) =
(R_{0,0} h,h)_H
= \| h \|^2_H, $$
for all $\delta\in \mathfrak{B}(\mathbb{R}^2)$, $h\in H$.  
Consequently, $\mu(\delta; h,g)$ is a sesquilinear (bilinear) functional with the norm less or equal to $1$. 
In fact, we may apply Theorem from~\cite[p. 64]{cit_300_A_G} (the proof of this theorem is valid for finite-dimensional
Hilbert spaces which are not ranked as Hilbert spaces in~\cite{cit_300_A_G}).
Therefore $\mu(\delta; h,g)$
admits the following representation:
\begin{equation}
\label{f3_45}
\mu(\delta; h,g) = (E(\delta) h,g)_H,\quad \delta\in \mathfrak{B}(\mathbb{R}^2),\ h,g\in H,
\end{equation}
where $E(\delta)$ is a linear bounded operator on $H$: $\| E(\delta) \| \leq 1$.
Observe that 
$$ (E(\delta) h,h)_H = \mu(\delta; h,h) \geq 0,\quad h\in H,\ \delta\in \mathfrak{B}(\mathbb{R}^2). $$
Therefore $E(\delta)\geq 0$, for all $\delta\in \mathfrak{B}(\mathbb{R}^2)$.
Thus, we have
\begin{equation}
\label{f3_47}
0\leq E(\delta)\leq E_H,\quad \delta\in \mathfrak{B}(\mathbb{R}^2).
\end{equation}
Notice that
$$ (E(\emptyset) h,g)_H = \mu(\emptyset; h,g) = 0, $$
$$ ( E( (0,2\pi]^2 ) h,g )_H = \mu((0,2\pi]^2; h,g) = \mu(\mathbb{R}^2; h,g) = (R_{0,0} h,g)_H = $$
$$ = (h,g)_H,\qquad h,g\in H. $$
Therefore
\begin{equation}
\label{f3_49}
E(\emptyset) = 0,\quad E((0,2\pi]^2) = E_H.
\end{equation}
For arbitrary $\delta_1,\delta_2\in\mathfrak{B}(\mathbb{R}^2)$, $\delta_1\cap\delta_2 = \emptyset$, and $h,g\in H$, we may write:
$$ (E(\delta_1\cup\delta_2) h,g)_H = \mu(\delta_1\cup\delta_2;h,g) = \mu(\delta_1;h,g) + \mu(\delta_2;h,g) = $$
$$ = (E(\delta_1) h,g)_H + (E(\delta_2) h,g)_H = ((E(\delta_1) + E(\delta_2)) h,g)_H, $$
and therefore
\begin{equation}
\label{f3_50}
E(\delta_1\cup\delta_2) = E(\delta_1) + E(\delta_2),\quad \delta_1,\delta_2\in \mathfrak{B}(\mathbb{R}^2):\ \delta_1\cap\delta_2 = \emptyset.
\end{equation}
Denote $K = \{ \delta\in \mathfrak{B}(\mathbb{R}^2):\ \delta\subseteq (0,2\pi]^2 \}$. 
By Neumark's theorem~\cite[p. 499]{cit_7500_R_S_N} we conclude that there exists 
a family $\{ F(\delta) \}_{\delta\in K}$ of operators of the orthogonal projection
in a Hilbert space $\widetilde H\supseteq H$ such that 
\begin{equation}
\label{f3_52}
F(\emptyset) = 0,\quad F((0,2\pi]^2) = E_{\widetilde H};
\end{equation}
\begin{equation}
\label{f3_55}
F(\delta_1\cap\delta_2) = F(\delta_1) F(\delta_2),\qquad  \delta_1,\delta_2\in K;
\end{equation}
\begin{equation}
\label{f3_57}
F(\delta\cup\widehat\delta) = F(\delta) + F(\widehat\delta),\qquad  \delta,\widehat\delta\in K:\ \delta\cap\widehat\delta=\emptyset;
\end{equation}
\begin{equation}
\label{f3_59}
E(\delta) = P^{\widetilde H}_H F(\delta) |_H,\qquad \delta\in K.
\end{equation}
Moreover, elements of the form $F(\delta) h$, $h\in H$, $\delta\in K$ determine $\widetilde H$.

Since $\mu$ is $\sigma$-additive, then by the latter property of $F$ we conclude that $F$ is weakly $\sigma$-additive.
In fact, let $\delta = \cup_{k=1}^\infty \delta_k$, where $\delta,\delta_k\in K$ and $\delta_i\cap\delta_j=\emptyset$,
$i,j\in\mathbb{N}: i\not= j$. For arbitrary $h,u\in H$ and $\widetilde\delta,\widehat\delta\in K$ we may write:
$$ \left(
\sum_{k=1}^N F(\delta_k) F(\widehat\delta) h, F(\widetilde\delta) u
\right)_{\widetilde H} =
\left(
\sum_{k=1}^N F(\delta_k\cap \widetilde\delta\cap \widehat\delta) h, u
\right)_{\widetilde H} = $$
$$ = \sum_{k=1}^N \left(
E(\delta_k\cap \widetilde\delta\cap \widehat\delta) h, u
\right)_{H} = 
\sum_{k=1}^N \mu\left( 
\delta_k\cap \widetilde\delta\cap \widehat\delta; h,u
\right) \rightarrow_{N\rightarrow +\infty}
$$
$$ \rightarrow_{N\rightarrow +\infty} 
\mu\left( 
\widetilde\delta\cap \widehat\delta \cap \left( \cup_{k=1}^{\infty} \delta_k \right) ; h,u
\right)
= 
\left( 
E\left( \widetilde\delta\cap \widehat\delta \cap \left( \cup_{k=1}^{\infty} \delta_k \right) \right) h,u
\right)_H = $$
$$ = \left( 
F\left( \widetilde\delta\cap \widehat\delta \cap \left( \cup_{k=1}^{\infty} \delta_k \right) \right) h,u
\right)_{\widetilde H} =
\left( 
F\left( \cup_{k=1}^{\infty} \delta_k \right) F(\widehat\delta) h, F(\widetilde\delta) u
\right)_{\widetilde H}.
$$
By the linearity we conclude that
$$ (S_N x,y)_{\widetilde H} \rightarrow_{N\rightarrow\infty} (S x,y)_{\widetilde H},\qquad x,y\in L, $$
where
$S_N := \sum_{k=1}^N F(\delta_k) = F\left(\cup_{k=1}^N \delta_k \right)$, $S := F\left( \cup_{k=1}^{\infty} \delta_k \right) = F(\delta)$,
$L := \mathop{\rm Lin}\nolimits \{ F(\delta) h:\ h\in H, \delta\in K \}$.
Choose arbitrary elements $h,g\in\widetilde H$. Since $L$ is dense in $\widetilde H$, there exist elements $h_k,g_k\in L$ such that
$\| h-h_k \| <\frac{1}{k}$, $\| g-g_k \| <\frac{1}{k}$, for all $k\in\mathbb{N}$. Observe that
$$ \left|
((S_N - S) h,g)_{\widetilde H} - ((S_N - S) h_k,g_k)_{\widetilde H}
\right| =
\left|
((S_N - S) h,g-g_k)_{\widetilde H} + 
\right. $$
$$ \left.
((S_N - S) (h-h_k),g_k)_{\widetilde H}
\right| \leq $$
$$\leq 2 \| h \| \| g-g_k \| + 2 \| h-h_k \| (\| g_k - g \| + \| g \|) \rightarrow_{k\rightarrow\infty} 0,\quad
(N\in\mathbb{N}). $$
For arbitrary $\varepsilon >0$ we may choose $k\in\mathbb{N}$ such that
$$ \left|
((S_N - S) h,g)_{\widetilde H} - ((S_N - S) h_k,g_k)_{\widetilde H}
\right| < \frac{\varepsilon}{2}. $$
There exists $\widehat N\in\mathbb{N}$ such that $N > \widehat N$ implies
$$ \left|
((S_N - S) h_k,g_k)_{\widetilde H} 
\right| < \frac{\varepsilon}{2}. $$
Then $\left|
((S_N - S) h,g)_{\widetilde H} \right| < \varepsilon$.
Therefore
\begin{equation}
\label{f3_60}
(S_N h,g)_{\widetilde H} \rightarrow_{N\rightarrow\infty} (S h,g)_{\widetilde H},\qquad h,g\in \widetilde H.
\end{equation}

\noindent
Define the following operator-valued functions:
\begin{equation}
\label{f3_63}
F_{1,t} = F((0,t]\times (0,2\pi]),\quad F_{2,t} = F((0,2\pi]\times (0,t]),\qquad t\in[0,2\pi].
\end{equation}
For $t<0$ we set $F_{1,t} = F_{2,t} = 0$, while for $t>2\pi$ we set $F_{1,t} = F_{2,t} = E_{\widetilde H}$.
Let us check that $\{ F_{j,t} \}$ is a spectral family on $[0,2\pi]$ such that $F_{j,0} = 0$; $j=1,2$.
By~(\ref{f3_52}) we see that
$F_{j,0} = 0$, $F_{j,2\pi} = E_{\widetilde H}$, $j=1,2$.
If $\lambda\leq\mu$, by~(\ref{f3_55}) we may write
$$ F_{1,\lambda} F_{1,\mu} = F((0,\lambda]\times (0,2\pi])  F((0,\mu]\times (0,2\pi]) = F((0,\lambda]\times (0,2\pi]) = $$
$$ = F_{1,\lambda}, $$
$$ F_{2,\lambda} F_{2,\mu} = F((0,2\pi]\times (0,\lambda]) F((0,2\pi]\times (0,\mu]) = F((0,2\pi]\times (0,\lambda]) = $$
$$ = F_{2,\lambda}. $$
It remains to check that $F_{j,t}$ is right-continuous ($j=1,2$).
For points $t\in (-\infty,0) \cup [2\pi,+\infty)$ it is obvious.
For arbitrary $t\in [0,2\pi)$; $t_k\in [0,2\pi):\ t_k > t$, $k\in\mathbb{N}$; $\{ t_k \}_1^\infty$ is decreasing and 
$t_k\rightarrow t$ as $k\rightarrow\infty$; 
and arbitrary $h,g\in\widetilde H$ we may write:
$$ ((F_{1, t_k}-F_{1,t}) h,g)_{\widetilde H} = 
\left(
F((t,t_k]\times (0,2\pi]) h,g
\right)_{\widetilde H}
= $$
$$ =
\left(
F\left( \cup_{n=1}^\infty ( (t_{n+1},t_n]\times (0,2\pi] ) \right) h,g
\right)_{\widetilde H} 
- $$
\begin{equation}
\label{f3_65}
-
\left(
F\left( \cup_{n=1}^{k-1} ( (t_{n+1},t_n]\times (0,2\pi] ) \right) h,g
\right)_{\widetilde H} 
\rightarrow_{k\rightarrow\infty} 
0. 
\end{equation}
Here we used the weak $\sigma$-additivity of $F$.
The monotone sequence of projections $\{ F_{1,t_n} \}_{n=1}^\infty$ converges in the strong operator topology to a
bounded operator. By~(\ref{f3_65}) we conclude that this operator is $F_{1,t}$.
If we would have $\lim_{u\to t+0} F_{1,u} h \not= F_{1,t} h$ for an element $h\in H$, then we could easily construct
a sequence $\{ t_k \}_{k=1}^\infty$ with above properties and satisfying
$\| F_{1,t_k} h - F_{1,t} h \| > \varepsilon$ with some $\varepsilon > 0$. This contradiction shows that
$F_{1,t}$ is right-continuous.
For $F_{2, t}$ we may use similar arguments.

By~(\ref{f3_55}) we may write
$$ F_{1,u} F_{2,v} = F((0,u]\times (0,2\pi]) F((0,2\pi]\times(0,v]) = F((0,u]\times (0,v]) = $$
\begin{equation}
\label{f3_69}
= F((0,2\pi]\times(0,v]) F((0,u]\times (0,2\pi]) = F_{2,v} F_{1,u},\qquad u,v\in [0,2\pi]. 
\end{equation}
Thus, $F_{1,u}$ and $F_{2,v}$ commute for all $u,v\in\mathbb{R}$.   
Set
\begin{equation}
\label{f3_70}
U_k = \int_0^{2\pi} e^{it} dF_{k,t},\qquad  k=1,2.
\end{equation}
Observe that $U_1$, $U_2$ are commuting unitary operators in $\widetilde H$.
By~(\ref{f3_45}),(\ref{f3_59}),(\ref{f3_69}) we may write
$$
\mu((a,b]\times(c,d]; h,h) = (E((a,b]\times(c,d]) h,h)_H =
(F((a,b]\times(c,d]) h,h)_{\widetilde H} = $$
\begin{equation}
\label{f3_75}
=
((F_{1,b} - F_{1,a}) (F_{2,d} - F_{2,c}) h, h)_{\widetilde H},\quad a,b,c,d\in [0,2\pi]:\ a<b,\ c<d,\ h\in H.
\end{equation}
By~(\ref{f3_20}) and (\ref{f3_75}) we conclude that
$$ \left(
\left. 
P^{\widetilde H}_H (E_{\widetilde H} + z_1 U_1) (E_{\widetilde H} - z_1 U_1)^{-1}
(E_{\widetilde H} + z_2 U_2) (E_{\widetilde H} - z_2 U_2)^{-1} \right|_H
h,h
\right)_H 
= $$
$$ 
=
\int_{\mathbb{R}^2}
\left(
\frac{1+z_1 e^{it_1}}{1-z_1 e^{it_1}}
\right)
\left(
\frac{1+z_2 e^{it_2}}{1-z_2 e^{it_2}}
\right)
d(F_{1,t_1} F_{2,t_2} h, h)_{\widetilde H} =  $$
$$ =
\int_{\mathbb{R}^2}
\left(
\frac{1+z_1 e^{it_1}}{1-z_1 e^{it_1}}
\right)
\left(
\frac{1+z_2 e^{it_2}}{1-z_2 e^{it_2}}
\right)
d\mu(\delta; h,h) =  ( R_{z_1,z_2} h, h )_H, $$
\begin{equation}
\label{f3_78}
z_1,z_2\in \mathbb{T}_e,\ h\in H.
\end{equation}
Consequently, $R_{z_1,z_2}$ is a generalized resolvent of a pair of isometric operators
$V_1 = V_2 = o_{H}$. Here $D(o_{H}) = \{ 0 \}$, $o_{H} 0 = 0$.
$\Box$

\begin{prop}
\label{p3_1}
Let an operator-valued function  $R_{z_1,z_2}$ be given, which depends on complex parameters
$z_1,z_2\in \mathbb{T}_e$ and which values are linear bounded operators defined on a (whole) Hilbert space $H$.
Let $V_1,V_2$ be closed isometric operators in $H$ which satisfy relation~(\ref{f1_1}).
Suppose that conditions 1)-3) of Theorem~\ref{t3_1} are satisfied. 
Suppose that conditions 1)-5) of Theorem~\ref{t1_3_p1_1} are satisfied with the choices
$V=V_1$, $R_\zeta = \frac{1}{2} \left( E_H + R_{\zeta,0} \right)$, and
$V=V_2$, $R_\zeta = \frac{1}{2} \left( E_H + R_{0,\zeta} \right)$.
Then $R_{z_1,z_2}$ is a generalized resolvent of a pair of isometric operators $V_1,V_2$.
\end{prop}
\textbf{Proof.}
Since all conditions of Theorem~\ref{t3_1} are satisfied, we can use the constructions from its proof. 
Thus, there exist commuting unitary operators $U_1,U_2$ in a Hilbert space $\widetilde H\supseteq H$ such that
\begin{equation}
\label{f3_80}
R_{z_1,z_2} =
\left. 
P^{\widetilde H}_H (E_{\widetilde H} + z_1 U_1) (E_{\widetilde H} - z_1 U_1)^{-1}
(E_{\widetilde H} + z_2 U_2) (E_{\widetilde H} - z_2 U_2)^{-1} \right|_H, 
\end{equation}
for $z_1,z_2\in \mathbb{T}_e$.
Then                    
\begin{equation}
\label{f3_87}
\frac{1}{2} \left( E_H + R_{\zeta,0} \right) =
\left. P^{\widetilde H}_H (E_{\widetilde H} - \zeta U_1)^{-1} \right|_H,\quad
\zeta\in \mathbb{T}_e;
\end{equation}
\begin{equation}
\label{f3_89}
\frac{1}{2} \left( E_H + R_{0,\zeta} \right) =
\left. P^{\widetilde H}_H (E_{\widetilde H} - \zeta U_2)^{-1} \right|_H,\quad
\zeta\in \mathbb{T}_e;
\end{equation}
and                        
\begin{equation}
\label{f3_93}
\left(
\frac{1}{2} \left( E_H + R_{\zeta,0} \right) h, g 
\right)_H =
\int_{\mathbb{R}} \frac{1}{1-\zeta e^{it}}
d(F_{1,t} h,g)_{\widetilde H},\quad
\zeta\in \mathbb{T}_e,\ h,g\in H;
\end{equation}
\begin{equation}
\label{f3_95}
\left(
\frac{1}{2} \left( E_H + R_{0,\zeta} \right) h, g 
\right)_H =
\int_{\mathbb{R}} \frac{1}{1-\zeta e^{it}}
d(F_{2,t} h,g)_{\widetilde H},\quad
\zeta\in \mathbb{T}_e,\ h,g\in H.
\end{equation}
Let us check that $U_1\supseteq V_1$.
Since conditions 1)-5) of Theorem~\ref{t1_3_p1_1} are satisfied with the choice
$V=V_1$, $R_\zeta = \frac{1}{2} \left( E_H + R_{\zeta,0} \right)$,
then choosing an arbitrary
$\zeta_0\in\mathbb{D}\backslash\{ 0 \}$ and
$L := (E_H - \zeta_0 V) D(V)$, we conclude that
conditions 1)-5) of Theorem~\ref{t1_2_p1_1} are satisfied, see the proof of Theorem~2 in~\cite{cit_750_Ch}.
Thus, $R_\zeta$ is a generalized resolvent of a closed isometric operator in a Hilbert space $H$ and
therefore $R_{\zeta_0}^{-1}$ exists and is a bounded operator on $H$.
Moreover, we have $D(V) = R_{\zeta_0} L$ (see the last formula on page~887 in~\cite{cit_750_Ch}).
By condition~1) of Theorem~\ref{t1_3_p1_1}
we have $V g = \frac{1}{\zeta_0} \left(
E_H - R_{\zeta_0}^{-1}
\right) g$, $g\in D(V)$.

\noindent
Thus, we can apply constructions from the proof of Theorem~1 in~\cite[p. 880]{cit_750_Ch}.
Notice that the above operator $V (=V_1)$ coincides with the operator $U$ defined by (30) in~\cite{cit_750_Ch}.
By formula~(26) in~\cite{cit_750_Ch} we may write:
\begin{equation}
\label{f3_97}
\left(
\frac{1}{2} \left( E_H + R_{\zeta,0} \right) h, g 
\right)_H =
\int_0^{2\pi} \frac{1}{1-\zeta e^{it}}
d (E_t h,g)_H,\quad
\zeta\in \mathbb{T}_e,\ h,g\in H.
\end{equation}
Comparing relations~(\ref{f3_93}) and~(\ref{f3_97}) we conclude that
\begin{equation}
\label{f3_102}
\int_0^{2\pi} \frac{1}{1-\zeta e^{it}}
d \left(
(E_t h,g)_H
-
(F_{1,t} h,g)_{\widetilde H}
\right) = 0,\quad
\zeta\in \mathbb{T}_e,\ h,g\in H.
\end{equation}
Therefore (see considerations on page~882 in~\cite[p. 883]{cit_750_Ch})
\begin{equation}
\label{f3_105}
\int_0^{2\pi} e^{it} d(E_t h,g)_H = \int_0^{2\pi} e^{it} d(F_{1,t} h,g)_{\widetilde H},\qquad  h,g\in H. 
\end{equation}
Then (cf.~\cite[p. 886]{cit_750_Ch})
$$ (V h,g)_H = \int_0^{2\pi} e^{it} d(E_t h,g)_H = \int_0^{2\pi} e^{it} d(F_{1,t} h,g)_{\widetilde H} = (U_1 h,g)_{\widetilde H}, $$  
\begin{equation}
\label{f3_107}
h\in D(V),\ g\in H. 
\end{equation}
Therefore $V h = P^{\widetilde H}_H U_1 h$, $h\in D(V)$. By $\| V h \| = \| U_1 h \|$ we get $U_1\supseteq V$.
Relation $U_2\supseteq V_2$ can be checked in the same manner.
By~(\ref{f3_80}) we see that
$R_{z_1,z_2}$ is a generalized resolvent of a pair $V_1,V_2$.
$\Box$

\begin{thm}
\label{t3_2}
Let an operator-valued function  $R_{z_1,z_2}$ be given, which depends on complex parameters
$z_1,z_2\in \mathbb{T}_e$ and which values are linear bounded operators defined on a (whole) Hilbert space $H$.
Let $V_1,V_2$ be closed isometric operators in $H$ which satisfy relation~(\ref{f1_1}).
$R_{z_1,z_2}$ is a generalized resolvent of a pair of
isometric operators $V_1,V_2$ if an only if the following conditions are satisfied:

\begin{itemize}

\item[1)] $R_{0,0} = E_H$;

\item[2)] $R_{z_1,z_2}^* = R_{ \frac{1}{\overline{z_1}}, \frac{1}{\overline{z_2}} }$, $z_1,z_2\in \mathbb{T}_e\backslash\{ 0 \}$;

\item[3)] For all $h\in H$, for the function $f(z_1,z_2) := (R_{z_1,z_2} h,h)_H$, $z_1,z_2\in\mathbb{T}_e$, there exist limits:
$$ f(\infty, z_2) := \lim_{z_1\to\infty} f(z_1,z_2),\quad f(z_1,\infty) := \lim_{z_2\to\infty} f(z_1,z_2),\quad z_1,z_2\in \mathbb{T}_e; $$
$$ f(\infty,\infty) = \lim_{z_2\to\infty} \lim_{z_1\to\infty} f(z_1,z_2), $$
and the extended by these relations function $f(z_1,z_2)$, $z_1,z_2\in \mathbb{T}_e\cup\{ \infty \}$ belongs to $H_2$.

\item[4)] $\frac{1}{2} \left( E_H + R_{\zeta,0} \right) (E_H - \zeta V_1) g = g$, for all $\zeta\in\mathbb{T}_e$, $g\in D(V_1)$;

\item[5)] $\frac{1}{2} \left( E_H + R_{0,\zeta} \right) (E_H - \zeta V_2) g = g$, for all $\zeta\in\mathbb{T}_e$, $g\in D(V_2)$.

\end{itemize}
\end{thm}
\textbf{Proof.}
\textit{Necessity}. The necessity of conditions~1)-3) follows from Theorem~\ref{t3_1}.
Repeating the arguments from the beginning of the proof of Proposition~\ref{p3_1} we conclude that
relations~(\ref{f3_87}), (\ref{f3_89}) hold. By condition~1) of Theorem~\ref{t1_3_p1_1} with 
$V=V_1$, $R_\zeta = \frac{1}{2} \left( E_H + R_{\zeta,0} \right)$, and
$V=V_2$, $R_\zeta = \frac{1}{2} \left( E_H + R_{0,\zeta} \right)$ it follows the validity of
conditions~4), 5) of the present theorem, respectively.

\noindent
\textit{Sufficiency.}
In order to apply Proposition~\ref{p3_1} it is sufficient to check that conditions 1)-5) of
Theorem~\ref{t1_3_p1_1} for the choices
$V=V_1$, $R_\zeta = \frac{1}{2} \left( E_H + R_{\zeta,0} \right)$, and
$V=V_2$, $R_\zeta = \frac{1}{2} \left( E_H + R_{0,\zeta} \right)$
are satisfied.
Condition~1) of Theorem~\ref{t1_3_p1_1} for these choices coincides with conditions 4),5) of the present theorem.
By Theorem~\ref{t3_1} and considerations in its proof $R_{z_1,z_2}$ is a generalized resolvent of $V_1=V_2=o_H$.
Then relations~(\ref{f3_80}), (\ref{f3_87}), (\ref{f3_89}) hold.
By Theorem~\ref{t1_3_p1_1} for $V=o_H$ and the above-mentioned choices of $R_\zeta$ we obtain that conditions~3),4),5) of
Theorem~\ref{t1_3_p1_1} are satisfied and they do not depend on $V$.
The required condition~2) of Theorem~\ref{t1_3_p1_1} for $V=V_1$, $R_\zeta = \frac{1}{2} \left( E_H + R_{\zeta,0} \right)$, and
$V=V_2$, $R_\zeta = \frac{1}{2} \left( E_H + R_{0,\zeta} \right)$ follows directly from
condition~1) of the present theorem.
$\Box$

\section{The case of commuting isometric and unitary operators.}

In this section we shall show how Theorem~\ref{t3_2} allows to parametrize generalized resolvents in the case of
commuting isometric and unitary operators.

Let $V_1=V$ be a closed isometric operator in a Hilbert space $H$, and $V_2=U$ be a unitary operator in $H$. Suppose that
relation~(\ref{f1_1}) holds. In our case it takes the following form:
\begin{equation}
\label{f4_10}
V U h = U V h,\quad h\in (U^{-1} D(V))\cap D(V).
\end{equation}
Suppose that there exist a Hilbert space $\widetilde H\supseteq H$ and
commuting unitary operators $U_1,U_2$ in $\widetilde H$, such that $U_1\supseteq V$, $U_2\supseteq U$.
Consider the corresponding generalized resolvent of a pair $V,U$:
$$ \mathbf{R}_{z_1,z_2} =
\left. 
P^{\widetilde H}_H U_1(z_1) U_2(z_2) \right|_H = 
P^{\widetilde H}_H U_1(z_1) U(z_2) = 
$$
\begin{equation}
\label{f4_15}
= P^{\widetilde H}_H U_1(z_1) |_H U(z_2) =
( -E_H + 2 \mathbf{R}_{z_1}(V) ) U(z_2),\qquad z_1, z_2\in \mathbb{T}_e,
\end{equation}
where $\mathbf{R}_{z_1}(V)$ is a generalized resolvent of the closed isometric operator $V$,
which corresponds to the unitary extension $U_1$.
On the other hand, we may write:                  
$$ \mathbf{R}_{z_1,z_2} =
\left. 
P^{\widetilde H}_H U_2(z_2) U_1(z_1) \right|_H = 
\left. 
P^{\widetilde H}_H U_2(z_2) |_H P^{\widetilde H}_H U_1(z_1) \right|_H =
$$
\begin{equation}
\label{f4_17}
=
U(z_2)( -E_H + 2 \mathbf{R}_{z_1}(V) ),\qquad z_1, z_2\in \mathbb{T}_e.
\end{equation}
Comparing relations~(\ref{f4_15}),(\ref{f4_17}) and simplifying we obtain that
\begin{equation}
\label{f4_19}
\mathbf{R}_{z_1}(V) (E_{H} - z_2 U)^{-1} =
(E_{H} - z_2 U)^{-1} \mathbf{R}_{z_1}(V),\qquad z_1, z_2\in \mathbb{T}_e.
\end{equation}
Therefore
\begin{equation}
\label{f4_22}
U \mathbf{R}_{z_1}(V) =
\mathbf{R}_{z_1}(V) U,\qquad z_1\in \mathbb{T}_e.
\end{equation}                                                         
By Chumakin's formula~(\ref{f1_15}) we may write:
\begin{equation}
\label{f4_25}
\mathbf{R}_{z_1}(V) = \left[ E_H - z_1 (V\oplus \Phi_{z_1}) \right]^{-1},\qquad z_1\in \mathbb{D},
\end{equation}
where $\Phi_{z_1}\in \mathcal{S}(\mathbb{D};N_0(V),N_\infty(V))$.
By~(\ref{f4_22}) and~(\ref{f4_25}) we obtain that
\begin{equation}
\label{f4_30}
(V\oplus \Phi_{z_1}) U = U (V\oplus \Phi_{z_1}),\qquad z_1\in \mathbb{D}.
\end{equation}
Here the equality for the case $z_1 = 0$ follows by the analyticity of $\Phi_{z_1}$.

\begin{thm}
\label{t4_1}
Let $V$ be a closed isometric operator in a Hilbert space $H$, and $U$ be a unitary operator in $H$. Suppose that
relation~(\ref{f4_10}) holds. 
Let $\mathcal{S}_{V,U}(\mathbb{D};N_0(V),N_\infty(V))$ be a set of all functions from $\mathcal{S}(\mathbb{D};N_0(V),N_\infty(V))$
which satisfy relation~(\ref{f4_30}).
Then the following statements hold:

\begin{itemize}

\item[(i)] The set of all generalized resolvents of a pair $V,U$ is non-empty if and only if
$\mathcal{S}_{V,U}(\mathbb{D};N_0(V),N_\infty(V))\not= \emptyset$;

\item[(ii)]
Suppose that $\mathcal{S}_{V,U}(\mathbb{D};N_0(V),N_\infty(V))\not= \emptyset$. An arbitrary generalized resolvent
of a pair $V,U$ has the following form:
\begin{equation}
\label{f4_35}
\mathbf{R}_{z_1,z_2} =
( -E_H + 2 \left[ E_H - z_1 (V\oplus \Phi_{z_1}) \right]^{-1} ) U(z_2),\qquad z_1\in\mathbb{D},\ z_2\in \mathbb{T}_e,
\end{equation}
where $\Phi_{z_1}\in \mathcal{S}_{V,U}(\mathbb{D};N_0(V),N_\infty(V))$,
and
\begin{equation}
\label{f4_37}
\mathbf{R}_{z_1,z_2} = \mathbf{R}_{ \frac{1}{\overline{z_1}}, \frac{1}{\overline{z_2}} }^*,\quad 
z_1\in\mathbb{D}_e,\ z_2\in \mathbb{T}_e\backslash\{ 0 \}.
\end{equation}
On the other hand, an arbitrary function $\Phi_{z_1}\in \mathcal{S}_{V,U}(\mathbb{D};N_0(V),N_\infty(V))$ defines
by relations~(\ref{f4_35}),(\ref{f4_37}) a generalized resolvent of a pair $V,U$
(for $z_1\in\mathbb{D}_e$, $z_2=0$ we define $\mathbf{R}_{z_1,z_2}$ by the weak continuity:
$\mathbf{R}_{z_1,0} = w.-\lim\limits_{z_2\to 0} \mathbf{R}_{z_1,z_2}$). 
Moreover, for different operator-valued functions from $\mathcal{S}_{V,U}(\mathbb{D};N_0(V),N_\infty(V))$ there correspond
different generalized resolvents of a pair $V,U$.
\end{itemize}
\end{thm}
\textbf{Proof.}
$(i):$
If the set of all generalized resolvents of a pair $V,U$ is non-empty, then by our considerations before the present theorem
we see that $\mathcal{S}_{V,U}(\mathbb{D};N_0(V),N_\infty(V))\not= \emptyset$.

\noindent  
On the other hand, suppose that $\mathcal{S}_{V,U}(\mathbb{D};N_0(V),N_\infty(V))\not= \emptyset$.
Choose an arbitrary function $\Phi_{z_1}\in \mathcal{S}_{V,U}(\mathbb{D};N_0(V),N_\infty(V))$.
Define a function $\mathbf{R}_{z_1,z_2}$ for
$(z_1,z_2)\in (\mathbb{D}\times\mathbb{T}_e) \cup (\mathbb{D}_e\times (\mathbb{T}_e\backslash\{ 0 \}))$
by relations~(\ref{f4_35}),(\ref{f4_37}).
Let $\mathbf{R}_{z_1}(V)$ be the generalized resolvent of $V$ corresponding to $\Phi_{\zeta_1}$ by Chumakin's formula.
By~(\ref{f4_30}) we obtain that relation~(\ref{f4_22}) holds for $z_1\in\mathbb{D}$. Therefore~(\ref{f4_22})
holds for all $z_1\in\mathbb{T}_e$, since the generalized resolvent $\mathbf{R}_{\zeta}(V)$ has the following property (\cite{cit_750_Ch}):
\begin{equation}
\label{f4_39}
\mathbf{R}_{\zeta}^*(V) = E_H - \mathbf{R}_{ \frac{1}{ \overline{\zeta} } }(V),\qquad \zeta\in \mathbb{T}_e\backslash\{ 0 \}.
\end{equation}
Consequently, relation~(\ref{f4_19}) holds and we may write:
\begin{equation}
\label{f4_45}
( -E_H + 2 \mathbf{R}_{z_1}(V) ) U(z_2) =
U(z_2) ( -E_H + 2 \mathbf{R}_{z_1}(V) ),\quad z_1,z_2\in\mathbb{T}_e. 
\end{equation}
By~(\ref{f4_45}) and our definition of $\mathbf{R}_{z_1,z_2}$, for arbitrary 
$z_1\in\mathbb{D}_e$, $z_2\in \mathbb{T}_e\backslash\{ 0 \}$ we may write:
$$ \mathbf{R}_{z_1,z_2} = \mathbf{R}_{ \frac{1}{\overline{z_1}}, \frac{1}{\overline{z_2}} }^* =
\left( U\left( \frac{1}{\overline{z_2}} \right) \right)^* \left( -E_H + 2 \mathbf{R}_{ \frac{1}{\overline{z_1}} }(V) \right)^* = 
$$  
\begin{equation}
\label{f4_49}
= U(z_2)
\left(
-E_H + 2 \mathbf{R}_{ z_1 }(V)
\right)
=
\left(
-E_H + 2 \mathbf{R}_{ z_1 }(V)
\right)
U(z_2). 
\end{equation}
Thus, for all $(z_1,z_2)\in (\mathbb{D}\times\mathbb{T}_e) \cup (\mathbb{D}_e\times (\mathbb{T}_e\backslash\{ 0 \}))$
we have the following representation:
\begin{equation}
\label{f4_50}
\mathbf{R}_{z_1,z_2} =
\left(
-E_H + 2 \mathbf{R}_{ z_1 }(V)
\right)
U(z_2).
\end{equation}
For a fixed $z_1\in\mathbb{D}_e$ by analyticity of $U(z_2)$ the following limit exists:
\begin{equation}
\label{f4_51}
w.-\lim\limits_{z_2\to 0} \mathbf{R}_{z_1,z_2}  = 
\left(
-E_H + 2 \mathbf{R}_{ z_1 }(V)
\right)
U(0) =: \mathbf{R}_{z_1,0}. 
\end{equation}
By~(\ref{f4_50}),(\ref{f4_51}),(\ref{f4_45}) we see that
\begin{equation}
\label{f4_55}
\mathbf{R}_{z_1,z_2} =
\left(
-E_H + 2 \mathbf{R}_{ z_1 }(V)
\right)
U(z_2)
=
U(z_2)
\left(
-E_H + 2 \mathbf{R}_{ z_1 }(V)
\right),\ z_1,z_2\in\mathbb{T}_e. 
\end{equation}
Let us check that $\mathbf{R}_{z_1,z_2}$ is a generalized resolvent of a pair $V,U$ by~Theorem~\ref{t3_2}.
The assumptions of Theorem~\ref{t3_2} with $V_1=V$, $V_2=U$, $R_{z_1,z_2}=\mathbf{R}_{z_1,z_2}$ and $H$ are satisfied.
Condition~1) of Theorem~\ref{t3_2} is satisfied, as well. 
By~(\ref{f4_55}) for arbitrary 
$z_1,z_2\in \mathbb{T}_e\backslash\{ 0 \}$ we may write:
$$ \mathbf{R}_{z_1,z_2}^* = 
\left(
-E_H + 2 \mathbf{R}_{ z_1 }(V)
\right)^* \left( U(z_2) \right)^*
=
\left( -E_H + 2 \mathbf{R}_{ \frac{1}{\overline{z_1}} }(V) \right)
U\left( \frac{1}{\overline{z_2}} \right) = 
$$
$$ = \mathbf{R}_{ \frac{1}{\overline{z_1}}, \frac{1}{\overline{z_2}} }. $$
Thus, condition~2) of Theorem~\ref{t3_2} is satisfied. 
By~(\ref{f4_55}) we see that
$$ \frac{1}{2} \left(
E_H + \mathbf{R}_{\zeta,0}
\right)
= \mathbf{R}_{\zeta}(V),\quad 
\frac{1}{2} \left(
E_H + \mathbf{R}_{0,\zeta}
\right)
= (E_H - \zeta U)^{-1},\quad \zeta\in\mathbb{T}_e. $$
Therefore condition~5) of Theorem~\ref{t3_2} is trivial and condition~4) of Theorem~\ref{t3_2}
follows from the property~1) of Theorem~\ref{t1_3_p1_1}.

\noindent
It remains to check condition~3) of Theorem~\ref{t3_2}. 
Since $\mathbf{R}_{ \zeta }(V)$ is a generalized resolvent of $V$, then there exists a unitary operator $Q\supseteq V$ in a Hilbert space
$\mathbf{H}\supseteq H$ such that
$$ \mathbf{R}_{ \zeta }(V) = \left. P^{\mathbf{H}}_H (E_{\mathbf{H}} - \zeta Q)^{-1}\right|_H,\qquad \zeta\in\mathbb{T}_e. $$
Then
$$ -E_H + 2 \mathbf{R}_{ z_1 }(V) = 
\left. P^{\mathbf{H}}_H Q(z_1) \right|_H,\qquad z_1\in\mathbb{T}_e. $$
Representation~(\ref{f4_55}) takes the following form:
\begin{equation}
\label{f4_57}
\mathbf{R}_{z_1,z_2} =
\left(
\left. P^{\mathbf{H}}_H Q(z_1) \right|_H
\right)
U(z_2)
=
U(z_2)
\left(
\left. P^{\mathbf{H}}_H Q(z_1) \right|_H
\right),\quad z_1,z_2\in\mathbb{T}_e. 
\end{equation}
Choose an arbitrary element $h\in H$. Set
$f(z_1,z_2) := (\mathbf{R}_{z_1,z_2} h,h)_H$, $z_1,z_2\in\mathbb{T}_e$.
Then
\begin{equation}
\label{f4_59}
f(z_1,z_2) =
\left(
Q(z_1)
(U(z_2) h), h \right)_{\mathbf{H}} 
=
\left(
U(z_2)
\left(
P^{\mathbf{H}}_H (Q(z_1) h)
\right), 
h \right)_H, 
\end{equation}
where $z_1,z_2\in\mathbb{T}_e$.
Since operator-valued functions $Q(z)$ and $U(z)$ are analytic at $\infty$, we conclude that the limits
in condition~3) of Theorem~\ref{t3_2} exist. Moreover, the limit values $f(\infty,z_2)$, $f(z_1,\infty)$, $f(\infty,\infty)$ may
be calculated by the formal substitution of $\infty$ in representations in~(\ref{f4_59}) using
$U(\infty) := -E_H$, $Q(\infty) := -E_{\mathbf{H}}$.
Thus, we may use representation~(\ref{f4_59}) for all values $z_1,z_2\in\mathbb{T}_e\cup\{ \infty \}$.

\noindent
Let us check that $f(z_1,z_2)$ ($z_1,z_2\in\mathbb{T}_e\cup\{ \infty \}$) belongs to the class $H_2$.
Holomorphy of $f(z_1,z_2)$ at $(z_1,z_2)$, $z_1,z_2\in\mathbb{T}_e\cup\{ \infty \}$ follows from 
holomorphy of $Q(z)$ and $U(z)$ at all points $z\in\mathbb{T}_e\cup\{ \infty \}$
and Hartogs's theorem.
By~(\ref{f4_5}),(\ref{f4_59}) it follows that
condition~(a) in the definition of $H_2$ holds.
Condition~(c) in the definition of $H_2$ follows by relation~(\ref{f4_59}).

\noindent
Let us check condition~(b) in the definition of $H_2$. 
Denote $W(z_1) = \left. P^{\mathbf{H}}_H Q(z_1) \right|_H$, $z_1\in \mathbb{T}_e\cup\{ \infty \}$. By~(\ref{f4_57}) we see that
\begin{equation}
\label{f4_60}
W(z_1) U(z_2)
=
U(z_2) W(z_1),\quad z_1,z_2\in\mathbb{T}_e\cup\{ \infty \}, 
\end{equation}
where the equality for infinite values of $z_1$ or $z_2$ holds trivially.
By~(\ref{f4_59}) we obtain that
$f(z_1,z_2) = (U(z_2) W(z_1) h,h)_H$, $z_1,z_2\in\mathbb{T}_e\cup\{ \infty \}$.
Choose arbitrary $z_1,z_2\in\mathbb{D}$ and write (cf.~\cite[p. 531]{cit_4000_Koranyi})
$$ f(z_1,z_2) - f(\overline{z_1}^{-1}, z_2) - f(z_1, \overline{z_2}^{-1}) + f(\overline{z_1}^{-1},\overline{z_2}^{-1}) = $$
\begin{equation}
\label{f4_62}
= (
( U(z_2) - U(\overline{z_2}^{-1}) ) ( W(z_1) - W(\overline{z_1}^{-1}) ) h, h)_H. 
\end{equation}
By~(\ref{f4_7}) it follows that operators $W(z_1) - W(\overline{z_1}^{-1})$, $U(z_2) - U(\overline{z_2}^{-1})$ are non-negative
bounded operators on $H$ (for $z_1,z_2=0$ it is trivial).
By~(\ref{f4_60}) we see that operators $W(z_1) - W(\overline{z_1}^{-1})$ and $U(z_2) - U(\overline{z_2}^{-1})$
commute.
Since the product of commuting bounded non-negative operators is non-negative, by~(\ref{f4_62}) we conclude
that condition~(b) in the definition of $H_2$ holds.
Consequently, $f(z_1,z_2)\in H_2$ and all conditions of Theorem~\ref{t3_2} are satisfied. By Theorem~\ref{t3_2}
we obtain that $\mathbf{R}_{z_1,z_2}$ is a generalized resolvent of the pair $V,U$.

\noindent
$(ii):$
If $\mathcal{S}_{V,U}(\mathbb{D};N_0(V),N_\infty(V))\not= \emptyset$, then by property~$(i)$ we see that
the set of all generalized resolvents of a pair $V,U$ is non-empty. 
Choose an arbitrary generalized resolvent $\mathbf{R}_{z_1,z_2}$ of a pair $V,U$.
By our considerations before the present theorem
we obtain that for $\mathbf{R}_{z_1,z_2}$ relation~(\ref{f4_35}) holds.
Relation~(\ref{f4_37}) follows by property~2) in Theorem~\ref{t3_2}.

\noindent
Choose an arbitrary function $\Phi_{z_1}\in \mathcal{S}_{V,U}(\mathbb{D};N_0(V),N_\infty(V))$.
Repeating considerations in the proof of condition~$(i)$ we conclude that a function $\mathbf{R}_{z_1,z_2}$, defined by 
relations~(\ref{f4_35}),(\ref{f4_37}), is a generalized resolvent of a pair $V,U$.

\noindent
For different operator-valued functions $\Phi_{z_1}$, $\widetilde \Phi_{z_1}$ from $\mathcal{S}_{V,U}(\mathbb{D};N_0(V),N_\infty(V))$ there correspond
different generalized resolvents of a closed isometric operator $V$. Suppose that $\Phi_{z_1}$, $\widetilde \Phi_{z_1}$
generate the same generalized resolvent
$\mathbf{R}_{z_1,z_2}$ of a pair $V,U$. Writing relation~(\ref{f4_35}) with $\Phi_{z_1}$ or $\widetilde \Phi_{z_1}$ and
$z_2=0$ we obtain  a contradiction.
$\Box$

In conditions of Theorem~\ref{t4_1} we additionally suppose that
\begin{equation}
\label{f4_70}
U D(V) = D(V).
\end{equation}
In this case condition~(\ref{f4_10}) implies $VU=UV$. Condition~(\ref{f4_30}) is equivalent to
\begin{equation}
\label{f4_75}
\Phi_{z_1} U g = U \Phi_{z_1} g,\qquad g\in H\ominus D(V),\ z_1\in\mathbb{D}.
\end{equation}
Observe that the function $\Phi_{z_1} = 0$ belongs to $\mathcal{S}(\mathbb{D};N_0(V),N_\infty(V))$ and
satisfies~(\ref{f4_75}). Thus, $\Phi_{z_1}\in \mathcal{S}_{V,U}(\mathbb{D};N_0(V),N_\infty(V))$ and therefore
the set of generalized resolvents of $V,U$ is non-empty.

Additionally suppose that $H$ is separable and there exists a conjugation $J$ on $H$ such that
\begin{equation}
\label{f4_79}
U J = J U^{-1},\quad J D(V) = R(V).
\end{equation}
Then 
\begin{equation}
\label{f4_82}
J (H\ominus D(V)) = H\ominus R(V).
\end{equation}
Denote
$$ U_0 := U|_{H\ominus D(V)}. $$
By the Godi\v{c}-Lucenko theorem~(\cite{cit_2000_GL}) for the unitary operator $U_0$ there exists the following representation:
\begin{equation}
\label{f4_83}
U_0 = KL, 
\end{equation}
where $K,L$ are two conjugations on a Hilbert space $H\ominus D(V)$.
Set
$$ \Theta = JK:\ H\ominus D(V)\rightarrow H\ominus R(V). $$
The operator $\Theta$ maps $H\ominus D(V)$ on the whole $H\ominus R(V)$ and $\Theta^{-1} = KJ$. 
By~(\ref{f4_79}),(\ref{f4_83}) we obtain that
\begin{equation}
\label{f4_85}
\Theta U g = U \Theta g,\qquad g\in H\ominus D(V).
\end{equation}
Then
\begin{equation}
\label{f4_86}
U \Theta^{-1} f = \Theta^{-1} U f,\qquad f\in H\ominus R(V).
\end{equation}

Let $\Psi_{z_1}$ be an arbitrary function from $\mathcal{S}(\mathbb{D};H\ominus D(V),H\ominus D(V))$
such that
\begin{equation}
\label{f4_87}
\Psi_{z_1} U_0 = U_0 \Psi_{z_1},\qquad z_1\in\mathbb{D}.
\end{equation}
Set
\begin{equation}
\label{f4_89}
\Phi_{z_1} = \Theta \Psi_{z_1},\qquad z_1\in\mathbb{D}.
\end{equation}
Observe that $\Phi_{z_1}$ belongs to $\mathcal{S}(\mathbb{D};N_0(V),N_\infty(V))$.
By~(\ref{f4_85}),(\ref{f4_87}) for arbitrary $g\in H\ominus D(V)$ and $z_1\in\mathbb{D}$ we may write:
$$ U \Phi_{z_1} g = U (\Theta (\Psi_{z_1} g )) = \Theta (U_0 (\Psi_{z_1} g ))  =  \Theta (\Psi_{z_1} ( U_0 g )) = \Phi_{z_1} U g. $$
Thus, $\Phi_{z_1}$ satisfies relation~(\ref{f4_75}).
Therefore $\Phi_{z_1}\in\mathcal{S}_{V,U}(\mathbb{D};N_0(V),N_\infty(V))$.

On the other hand, choose an arbitrary $\Phi_{z_1}\in\mathcal{S}_{V,U}(\mathbb{D};N_0(V),N_\infty(V))$.
Then $\Phi_{z_1}$ belongs to $\mathcal{S}(\mathbb{D};N_0(V),N_\infty(V))$ and satisfies relation~(\ref{f4_75}).
Set
\begin{equation}
\label{f4_91}
\Psi_{z_1} = \Theta^{-1} \Phi_{z_1},\qquad z_1\in\mathbb{D}.
\end{equation}
Then relation~(\ref{f4_89}) holds. Observe that $\Psi_{z_1}\in\mathcal{S}(\mathbb{D};H\ominus D(V),H\ominus D(V))$.
Fix an arbitrary $z_1\in\mathbb{D}$. By~(\ref{f4_75}),(\ref{f4_86}) for arbitrary $g\in H\ominus D(V)$ we may write:
$$ \Psi_{z_1} U_0 g = \Theta^{-1} (\Phi_{z_1} (U g)) = \Theta^{-1} (U (\Phi_{z_1} g));  $$
$$ U_0 \Psi_{z_1} g =  U (\Theta^{-1} (\Phi_{z_1} g)) = \Theta^{-1} (U (\Phi_{z_1} g)). $$
Therefore relation~(\ref{f4_87}) holds.


\begin{center}
{\large\bf 
Characteristic properties for a generalized resolvent of a pair of commuting isometric operators.}
\end{center}
\begin{center}
{\bf S.M. Zagorodnyuk}
\end{center}

In this paper we consider a notion of a generalized resolvent for a pair of commuting isometric operators in a Hilbert space
$H$. Characteristic properties of the generalized resolvent are obtained. 
}

\end{document}